\documentclass[12pt]{article}

\usepackage{amsmath}
\usepackage{amsfonts}
\usepackage{amssymb}
\usepackage{graphicx, color}
\usepackage{hyperref}
\usepackage{epsfig}
\usepackage{epstopdf} 

\definecolor{darkgreen}{rgb}{0,0.55,0}

\newtheorem{proposition}{Proposition}[section]
\newtheorem{theorem}{Theorem}[section]
\newtheorem{lemma}[theorem]{Lemma}
\newtheorem{corollary}[theorem]{Corollary}
\newtheorem{remark}[theorem]{Remark}

\newtheorem{definition}{Definition}

\def\phi{{\varphi}}

\DeclareSymbolFont{AMSb}{U}{msb}{m}{n}
\DeclareMathSymbol{\N}{\mathbin}{AMSb}{"4E}
\DeclareMathSymbol{\Z}{\mathbin}{AMSb}{"5A}
\DeclareMathSymbol{\R}{\mathbin}{AMSb}{"52}
\DeclareMathSymbol{\Q}{\mathbin}{AMSb}{"51}
\DeclareMathSymbol{\I}{\mathbin}{AMSb}{"49}
\DeclareMathSymbol{\C}{\mathbin}{AMSb}{"43}

\DeclareMathOperator*{\esssup}{ess\,sup}

\newcommand{\calH}{{\mathcal H}}

\begin{document}
\title{Existence and structure of minimizers of  least gradient problems}

\author{{Amir Moradifam\footnote{Department of Mathematics, University of California, Riverside, CA, USA. E-mail: moradifam@math.ucr.edu.  }
}}
\date{\today}

\smallbreak \maketitle

\begin{abstract}
We study existence of minimizers of the general least gradient problem
\[\inf_{u \in BV_f} \int_{\Omega}\varphi(x,Du),\]
where $BV_f=\{u \in BV(\Omega): \ \ u|_{\partial \Omega}=f\}$, $f\in L^{1}(\partial \Omega)$, and $\varphi(x,\xi)$ is convex, continuous, and homogeneous function of degree $1$ with
respect to the $\xi$ variable. It is proven that there exists a divergence free vector field $T\in (L^{\infty}(\Omega))^n$ that determines the structure of level sets of all (possible) minimizers, i.e. $T$ determines $\frac{Du}{|Du|}$, $|Du|-$ a.e. in $\Omega$, for all minimizers $u$.  We  also prove that every minimizer of the above least gradient problem is also a minimizer of 
\[\inf_{u\in \mathcal{A}_f} \int_{\R^n}\varphi(x,Du),\]
where $\mathcal{A}_f=\{v\in BV(\R^n): \ \ v=f \ \ \hbox{on}\ \ \Omega^c\}$ and $f\in W^{1,1}(\R^n)$ is a compactly supported extension of $f\in L^1(\partial \Omega)$, and show that $T$ also determines the structure of level sets of all minimizers of the latter problem. This relationship between minimizers of the above two least gradient problems could be exploited to obtain information about existence and structure of minimizers of the former problem from that of the latter, which always exist.

\end{abstract}

\section{Introduction and Statement of Results}
Let $\Omega$ be a bounded open set in $\R^n$ with Lipschitz boundary and $\varphi: \Omega\times \R^n \rightarrow \R$ be a continuous function satisfying the following conditions:\\

($C_1$) There exists $\alpha>0$ such that $0 \leq \varphi(x,\xi) \leq \alpha |\xi|$ for all $x\in \Omega$ and $\xi \in \R^n$. \\

($C_2$) $\xi \mapsto \varphi(x,\xi)$ is a norm for every $x$.\\ \\
For any $u\in BV_{loc}(\R^n)$ let $\varphi(x,Du)$ denote the measure
defined by
\begin{equation}\label{varphi.def}
\int_A \varphi(x,Du) \ = \ \int_A \varphi(x, v^u(x)) |Du|
\qquad\mbox{ for any bounded Borel set $A$},
\end{equation}
where $|Du|$ is the total variation measure associated to the vector-valued
measure $Du$, and $v^u$ denotes the Radon-Nikodym derivative $v^{u}(x)=\frac{d\, Du}{d\, |Du|}$.
Standard  facts about $BV$ functions imply that (see \cite{AB}) if $U$ is an open set,
then
\begin{equation}\label{varphi.def1}
\int_{U} \varphi(x,Du) = \sup \{ \int_{U} u \nabla \cdot Y dx\ \ : \ \ Y \in C^{\infty}_{c} (U; \R^n), \ \ \sup \varphi^0(x, Y(x)) \leq 1 \}, 
\end{equation}
where $\varphi^0(x, \cdot)$ denotes the norm on $\R^n$ dual to $\varphi(x, \cdot)$, defined by
\[
\varphi^0(x,\xi) := \sup \{ \xi \cdot p : \varphi(x,p) \le 1\}.
\]
Since $\varphi$ satisfies ($C_1$), the dual $\varphi^0(x, \cdot)$ can be equivalently defined by 
\begin{equation}\label{varphi0}
\varphi^0(x,\xi) =\sup \{\frac{\xi \cdot p}{\varphi(x,p)}: \ \ p\in \R^n\}, 
\end{equation}
(see (2.17) in \cite{AB}). 
For $u\in BV(\Omega)$, $\int_{\Omega}\varphi(x, Du)$ is called the $\varphi$-total variation of $u$ in $\Omega$. 

In this paper, we study existence and structure of  minimizers of the general least gradient problem 
\begin{equation}\label{LTVProb}
\inf_{v\in BV_f(\Omega)} \int_{\Omega} \varphi(x, Dv)
\end{equation}
where
$f\in L^1(\partial \Omega)$ and 
\[
BV_f(\Omega):=\{v\in BV(\Omega): 
\ \ \hbox{for a.e. } x\in \partial \Omega, \ \ 
\lim_{r\to 0} \ \esssup_{y\in \Omega, |x-y|<r} |f(x) - v(y)| = 0 \ \}.
\]

Least gradient problems naturally arise in conductivity imaging. In \cite{HMN} the author and collaborators presented a method for recovering the conformal factor of an anisotropic conductivity matrix in a known conformal class from one interior measurement.  More precisely, assume that the matrix valued conductivity $\sigma(x)$ is of the form
\[\sigma(x) =c(x)\sigma_0(x)\]
where $c(x)\in C^{\alpha}(\Omega)$ is a positive scalar valued function and $\sigma_0 \in C^{\alpha}(\Omega, Mat(n, \R^n))$ is  a known positive definite symmetric matrix valued function. In medical imaging  $\sigma_0$ can be determined using Diffusion Tensor Magnetic Resonance Imaging. In \cite{HMN} the authors showed that the corresponding voltage potential $u$ is the unique solution of the least gradient problem 
\[ \hbox{argmin} \{ \int_{\Omega} \varphi(x, Dv): \ \ u \in {BV(\Omega)}, \ \ u|_{\partial \Omega}=f \},\]
where $\varphi$ is given by 
\begin{equation}\label{varphi}
\varphi(x,\xi)=a(x)\left( \sum_{i,j=1}^{n}\sigma_0^{ij}(x)\xi _i \xi_j\right)^{1/2},
\end{equation}
\begin{equation}\label{aJJ}
a=\sqrt{\sigma_0^{-1} J \cdot J},
\end{equation}
and $J$ is the current density vector field generated by imposing the voltage $f$ at $\partial \Omega$.  Once $u$ is determined the function $c(x)$ can easily be calculated. Recovering isotropic conductivities is a special case of the above formulation where $\sigma_0$ is the identity matrix and the weight $a$ is the magnitude of the induced current density vector field. See \cite{MNT, MNTa_SIAM, NTT07, NTT08, NTT10, NTT11} for applications of least gradient problems in imaging isotropic conductivities.

Any function $f\in L^1(\partial \Omega)$ can be extended to a compactly supported function in $W^{1,1}_c(\R^n)$ with inner and outer trace $f$ on $\partial \Omega$ (see for example \cite{G}). Throughout the paper we will denote this function by $f$ again and assume that $f\in L^1(\partial \Omega)$ is the restriction of a function $f\in W^{1,1}_c(\R^n)$. We will frequently switch between writing $f\in L^1(\partial \Omega)$ and $f\in W^{1,1}_c(\R^n)$ depending on the context. 

It is well-known that, the least gradient problem  \eqref{LTVProb} may not have a minimizer in $BV_f(\Omega)$(see \cite{sternberg_ziemer92}, \cite{JMN}, \cite{MRD}).  To see this, suppose $\{u_n\}_{n=1}^{\infty}$ is a minimizing sequence of (\ref{LTVProb}). Since $BV(\Omega)\hookrightarrow L^1_{loc}(\Omega)$, $I(v)= \int_{\Omega} \varphi(x, Dv)$ is coercive in $BV(\Omega)$ (a consequence of $C_1$) and weakly lower semicontinuous (see \cite{JMN} for more details), it follows from standard arguments that $\{u_n\}_{n=1}^{\infty}$ has a subsequence converging strongly in $L^1_{loc}$ to a function $\tilde{u}\in BV(\Omega)$ with 
\[\int_{\Omega}\varphi(x, D\tilde{u}) \leq \inf_{v\in BV_f(\Omega)} \int_{\Omega} \varphi(x, Dv).\] 
However, in general, the trace $\tilde{u} |_{\partial \Omega}$ on $\partial \Omega$ may not be equal to $f$, leading to possible nonexistence for the problem \eqref{LTVProb}. A natural question one may ask is whether it is possible to deduce information about existence, multiplicity, and structure of minimizers of (\ref{LTVProb}) in $BV_f$ from the knowledge of a limit $\tilde{u}\in BV(\Omega)$ of a minimizing sequence $\{u_n\}_{n=1}^{\infty}$, which may not have the trace $f$ on $\partial \Omega$.  One of the main objectives of the paper is to answer this question.  We shall show that $\tilde{u}$ reveals fundamental information about existence and the structure of level sets of the minimizers of (\ref{LTVProb}).

 Define 
\[\mathcal{A}_f:= \{v \in BV(\R^n): \ \ v=f \ \ \hbox{on}\ \ \Omega^c\},\]
and note that $BV_f \subsetneqq \mathcal{A}_f$ and $BV_f \hookrightarrow \mathcal{A}_f$ in the sense that any element $v$ of $BV_f(\Omega)$ is the restriction to $\Omega$ of a unique element of $\mathcal{A}_f$. It follows from the above argument that any minimizing sequence  $\{v_n\}_{n=1}^{\infty}$ of \eqref{LTVProb} has a subsequence converging strongly in $L^1_{loc}$ to a function  $w \in \mathcal{A}_f$ satisfying
\[\int_{\Omega}\varphi(x, Dw) \leq \inf_{v\in \mathcal{A}_f} \int_{\Omega} \varphi(x, Dv).\] 
Hence $w$ is a minimizer of the least gradient problem 
\begin{equation}\label{LTVProbRelaxed}
\inf_{v\in \mathcal{A}_f} \int_{\Omega} \varphi(x, Dv).
\end{equation}
One of our main goals is to study the relation between minimizers of \eqref{LTVProbRelaxed} (which always exist) and the existence of minimizers of \eqref{LTVProb}.  We shall first prove that any minimizer of \eqref{LTVProb} is also a minimizer of \eqref{LTVProbRelaxed}. 

\begin{proposition}\label{SameValue}
Let $\Omega \subset \R^n$ be a bounded open set with Lipschitz boundary and assume $\varphi: \Omega \times \R^n \rightarrow \R$ be a continuous function satisfying the condition $(C_1)$ and $(C_2)$, and $f\in W^{1,1}_c(\R^n)$.  Then 
\begin{eqnarray*}
\min _{v\in \mathcal{A}_f}  \left( \int_{\Omega} \varphi(x, Dv)+\int_{\partial \Omega}\varphi(x,\nu_{\Omega})|f-v|ds \right)=   \inf_{v\in BV_f(\Omega)} \int_{\Omega} \varphi(x, Dv). \\ \\
\end{eqnarray*}
\end{proposition}

Next we prove that all minimizers of the least gradient problems \eqref{LTVProbRelaxed} and \eqref{LTVProb} have the same level set structure, confirming an observation of  Maz\'{o}n, Rossi, and De Le\'{o}n \cite{MRD} (see Remark 2.8 in \cite{MRD}).

\begin{theorem}\label{Structue}
Let $\Omega \subset \R^n$ be a bounded open set with Lipschitz boundary and assume $\varphi: \Omega \times \R^n \rightarrow \R$ be a continuous function satisfying the condition $(C_1)$ and $(C_2)$, and $f\in W^{1,1}_c(\R^n)$.  Then there exists a divergence free vector field $T \in (L^{\infty}(\Omega))^n$ with $\varphi^0(x,T)\leq 1$ a.e. in $\Omega$ such that every minimizer $w$ of \eqref{LTVProb} or \eqref{LTVProbRelaxed} satisfies 
\begin{equation}\label{directionInside}
\varphi(x,\frac{Dw}{|Dw|})=T\cdot \frac{Dw}{|Dw|}, \ \ |Dw|-a.e. \ \ \hbox{in} \ \ \Omega
\end{equation}
and
\begin{equation} \label{Boundary1}
\varphi(x,\nu_{\Omega})= [T, \hbox{sign} (f-w) \nu_{\Omega}], \ \ \mathcal{H}^{n-1}-a.e. \ \ \hbox{on} \ \ \partial \Omega.\\ \\ \\ 
\end{equation}
\end{theorem}
\vspace{1cm}

The above theorem asserts that a fixed divergence free vector field $T$ determines the structure of the level sets of all minimizers of  the least gradient problems \eqref{LTVProb} and (\ref{LTVProbRelaxed}). More precisely, since $\varphi^0(x,T)\leq 1$ a.e. in $\Omega$, we have
\[\varphi(x,p)\geq T\cdot p  \]
for every $p\in S^{n-1}$ and a.e.  $x\in  \Omega$. Thus it follows from \eqref{directionInside} that $|Dw|$-a.e., $p=\frac{Dw}{|Dw|}$ maximizes 
\[\frac{T \cdot p}{\varphi(x,p)}\]
among all $p\in S^{n-1}$, determining  $\frac{Dw}{|Dw|}$, $|Du|$-a.e. in $\Omega$. Theorem \ref{Structue} should be compared to the results  in \cite{Mazon2}. 

On the other hand, the condition \eqref{Boundary1} determines the set of possible jumps of a minimizer $u$ on $\partial \Omega$. To see this, suppose the trace of $T$ can be represented by a function $T_{tr}\in ( L^{\infty}(\partial \Omega))^n$. Then \eqref{Boundary1} implies that, up to a set with  $\mathcal{H}^{n-1}$-measure zero, 

\[\{x\in \partial \Omega: w|_{\partial \Omega}> f\} \subseteq \{x\in \partial \Omega: \varphi(x,\nu_{\Omega}(x))=T_{tr}\cdot \nu_{\Omega}\},\]
and similarly 
\[\{x\in \partial \Omega: w|_{\partial \Omega}< f\} \subseteq \{x\in \partial \Omega: \varphi(x,\nu_{\Omega}(x))=-T_{tr}\cdot \nu_{\Omega}\},\]
for every minimizer $w$ of \eqref{LTVProbRelaxed}. The above conclusions are more explicit in the following corollary of Theorem \ref{Structue}. 

\begin{corollary}\label{Structure1}
Let $\Omega \subset \R^n$ be a bounded open set with Lipschitz boundary and assume that $a \in C(\overline{\Omega}) $ is a non-negative function, and $f\in W^{1,1}_c(\R^n)$.  Then there exists a divergence free vector field $T \in (L^{\infty}(\Omega))^n$ with $|T|\leq a$ a.e. in $\Omega$ such that  every minimizer $w\in \mathcal{A}_f$ of 
\begin{equation}\label{LGa}
\inf_{v\in \mathcal{A}_f} \int_{\Omega} a|Dv|,
\end{equation}
 satisfies 
\begin{equation}
T\cdot \frac{Dw}{|Dw|}=|T|=a, \ \ |Dw|-a.e. \ \ \hbox{in} \ \ \Omega,
\end{equation}
and
\begin{equation} \label{Boundary2}
a= [T, \hbox{sign} (f-u) \nu_{\Omega}], \ \ \mathcal{H}^{n-1}-a.e. \ \ \hbox{on} \ \ \partial \Omega.
\end{equation}
\end{corollary}
\vspace{1cm}

Corollary \ref{Structure1} asserts that there exists a divergence free vector field $T$ such that for every minimizer $u$ of \eqref{Boundary1} the vector field $\frac{Dw}{|Dw|}$ is parallel to $T$, $|Dw|$-a.e. in $\Omega$. See Section 5 in \cite{JMN} for an example of a least gradient problem that has infinitely many minimizers, all of which have the same level set structure. Moreover, if the trace of $T$ can be represented by a function $T_{tr}\in ( L^{\infty}(\partial \Omega))^n$, then up to a set with $\mathcal{H}^{n-1}$-measure zero 
\[\{x\in \partial \Omega: w|_{\partial \Omega}> f\} \subseteq \{x\in \partial \Omega: T_{tr} \cdot \nu_{\Omega}=|T_{tr}|\},\]
and similarly 
\[\{x\in \partial \Omega: w|_{\partial \Omega}< f\} \subseteq \{x\in \partial \Omega:  T_{tr} \cdot \nu_{\Omega}=-|T_{tr}|\}.\]
In other words $w|_{\partial \Omega}=f$, $\mathcal{H}^{n-1}$-a.e. in 
\[\{ x\in \partial \Omega: |T_{tr}\cdot \nu_{\Omega}|<|T_{tr}|\},\]
for every minimizer $w$ of \eqref{LTVProbRelaxed}. \\

\begin{remark} \label{nonexistence}  Suppose that assumptions of Corollary \ref{Structure1} hold and let $w \in \mathcal{A}_f$ be a minimizer of \eqref{LGa} with $w|_{\Gamma}\neq f$, for some open subset $\Gamma$ of $\partial \Omega$. Also let $T$ be the vector field in the statement of Corollary \ref{Structure1}, and assume that $T$ is continuous in a neighborhood of $\Gamma$, i.e. $T \in C(\Omega \cap \mathcal{O}) \cup \Gamma)$, where $\mathcal{O}$ is an open set of $\R^n$ containing $\Gamma$. If $\tilde{w}$ is another minimizer of \eqref{LTVProbRelaxed} which is locally $C^1$ near $\Gamma$ and satisfies $\tilde{w}|_{\Gamma}=f$, then $f$ must be constant along $\Gamma$. Indeed, since $w$ has a jump on $\Gamma$, it follows from \eqref{Boundary2} that $T$ is parallel to $\nu_{\Omega}$ on $\Gamma$. Therefore, by Corollary \ref{Structure1}, $\nabla \tilde{w}$ is also parallel to $\nu_{\Omega}$ on $\Gamma$. Thus $\tilde{w}$ must be constant on the jump set $\Gamma \subset \partial \Omega$ of $w$. In particular, if $f$ is not constant on every open connected component of the jump set $\Gamma \subset \partial \Omega$ of $w \in \mathcal{A}_f$, then \eqref{LGa} does not have a minimizer in $BV_f$ that is locally $C^1$ near $\Gamma$. 
\end{remark}

In what follows we are concerned with sufficient conditions to guarantee that every minimizer $w\in \mathcal{A}_f$ of (\ref{LTVProbRelaxed}) belongs to $BV_f(\Omega)$ and therefore is also a minimizer of the least gradient problem \eqref{LTVProb}.  In \cite{JMN}, the author and collaborators showed if $f\in C(\partial \Omega)$ and $\partial \Omega$ satisfies following geometric hypothesis, then every minimizer of \eqref{LTVProbRelaxed} is also a minimizer of \eqref{LTVProb} (see Theorem 1.1. in \cite{JMN}).

For $u\in BV(\Omega)$, $\int_{\R^n}\varphi(x, Du)$ is called the $\varphi$-total variation of $u$ in $\R^n$. Also if  $E$ is a Borel subset of $\R^n$, then we shall write $P_{\varphi}(E; \R^n)$ to denote the $\varphi$-perimeter of $E$ in $R^n$, defined by
\begin{equation}\label{perimeter}
P_{\varphi}(E;\R^n):= \int_{\R^n} \varphi(x, D \chi_E) ,
\end{equation}
where $\chi_E$ is the characteristic function of $E$. Note that if $\partial E$ is smooth enough, then
\[
P_{\varphi}(E;\R^n):= \int_{\partial E }  \varphi(x, \nu_E(x)) \,d\calH^{n-1}\qquad
\mbox{$\nu_E :=$  outer unit normal,}
\]
which is
a generalized inhomogeneous, anisotropic area of $\partial E$ in $\R^n$.

If $V$ is a measurable subset of $\R^n$, we will write 
\begin{equation}
V^{(1)}:= \{x\in \R^n :
\lim_{r\to 0} \frac {\calH^n(B(r,x)\cap V)}{\calH^n(B(r))} = 1 \}. 
\end{equation}
\begin{definition} Let $\Omega \subset \R^n$ be a bounded Lipschitz domain and $\varphi: \Omega \times \R^n \rightarrow \R$ is a continuous function that satisfies $C_1$ and $C_2$. We say that $ \Omega$ satisfies the  barrier condition if for every $x_0\in \partial \Omega$ and $\epsilon>0$ sufficiently small, if $V$
minimizes $P_{\varphi}(\,\cdot\,; \R^n)$ in
\begin{equation}\label{SBC-Set}
\{W \subset \Omega : W \setminus  B(\epsilon, x_0)=\Omega \setminus B(\epsilon, x_0)\},
\end{equation}
then
\[\partial V^{(1)} \cap \partial \Omega \cap B(\epsilon, x_0)= \emptyset.\]
\label{def.barrier}\end{definition}

When $\varphi(x,\xi)= |\xi|$, the above condition is equivalent to those introduced by Sternberg and Ziemer (see (3.1) and (3.2) in \cite{sternberg_ziemer92}), at
least for smooth sets. 

\begin{remark}
In \cite{JMN}, it is proved that if $\varphi \in C^1$ and  $\partial \Omega$ is sufficiently smooth, then $\Omega$ satisfies the barrier condition provided
\begin{equation}
-\sum_{i=1}^n \partial_{x_i} \varphi_{\xi_i} (x, Dd(x)) >0 \quad\mbox{ on a dense subset of 
$\partial \Omega$},
\label{barrier.suff}\end{equation}
where 
\[
d(x) := \begin{cases}
\operatorname{dist}(x, \partial \Omega) &\mbox{ if }x\in \Omega\\
-\operatorname{dist}(x, \partial \Omega) &\mbox{ if not}.
\end{cases}
\]
\end{remark}
\vspace{.2cm}

\begin{theorem}\label{existence.Theo}
Suppose that $\varphi: \R^n \times \R^n \rightarrow \R$ is a continuous function that satisfies $C_1$ and $C_2$ in a bounded Lipschitz domain $\Omega\subset\R^n$.
If $\Omega$ satisfies the barrier condition  with respect to $\varphi$ and $f\in W^{1,1}_c(\R^n)$ is continuous at $\mathcal{H}^{n-1}$-a.e. $x \in \partial \Omega$, then every minimizer $w\in \mathcal{A}_f$ of (\ref{LTVProbRelaxed}) is also a minimizer of (\ref{LTVProb}). In particular, the least gradient problem (\ref{LTVProb}) has a minimizer in $BV_f(\Omega)$. 
\end{theorem}

The proof of the above theorem follows from a slight modification of the proof of Theorem 1.1 in \cite{JMN}, and will not be presented here. For the case $\varphi(x,\xi)=|\xi|$ and $f\in C(\partial \Omega)$, Theorem \ref{existence.Theo} reduces to the existence result of Sternberg and Ziemer in \cite{sternberg_ziemer92} which is the first result in this direction (see also \cite{sternberg_ziemer}, \cite{sternbergZiemer93}, and \cite{Gorny}).

 In \cite{ST}, Spradlin and Tamasan considered the case $\varphi(x,\xi)= |\xi|$ and presented an example of an $L^1$ function on the unite disk that satisfies the barrier condition but is not the trace of a function of least gradient. This function is the characteristic of a Cantor set which is discontinuous on a set of positive measure on the unit circle, and hence Theorem \ref{existence.Theo} does not apply. Indeed the example of Spradlin and Tamsan shows that Theorem \ref{existence.Theo} is sharp.

\section{Proofs} 
Let $\Omega$ be a bounded open set in $\R^n$ with Lipschitz boundary, and $f\in L^1(\partial \Omega)$ be the restriction of a compactly supported function (denoted by $f$ again) in $W^{1,1}(\R^n)$ with inner and outer trace $f$ on $\partial \Omega$. Define 
\[\mathcal{A}_f:= \{w \in BV(\R^n): \ \ w=f \ \ \hbox{on}\ \ \Omega^c\},\]
and note that $BV_f \hookrightarrow \mathcal{A}_f$ in the sense that any element $v$ of $BV_f(\Omega)$ is the restriction to $\Omega$ of a unique element of $\mathcal{A}_f$. The problem \eqref{LTVProb} may not have a solution, but as argued in the introduction (\ref{LTVProbRelaxed}) always has a solution. 

Let $E: (L^1(\Omega))^n\rightarrow \R$ and $G: W^{1,1}_0(\Omega) \rightarrow \R$ be defined as follows
\begin{equation}\label{FandG}
E(P):=\int_{\Omega} \varphi (x, P+\nabla f)dx, \ \ G(u)\equiv 0.
\end{equation}
Then the problem (\ref{LTVProbRelaxed}) can be written as 
\[(P) \ \ \ \ \ \inf_{u \in W^{1,1}_0(\Omega)} E(Du)+G(u).\]
By Fenchel duality (see Chapter III in \cite{ET}) the dual problem is given by  
\begin{equation*}
(P^*) \hspace{.5cm}\sup_{V \in (L^{\infty}(\Omega))^n} \{ -E^* (V)-G^*(-\nabla \cdot V)\}.
\end{equation*}
Recall that the Legendre-Fenchel transform  $E^*: (L^{\infty}(\Omega))^n \rightarrow \R$ is 
\[E^*(V)=\sup  \{\langle V, P \rangle -E(P): \ \  P \in  (L^1(\Omega))^n\}.\]
One can easily compute $G^*: W^{-1,\infty} (\Omega)\rightarrow \R$:
\begin{equation*}G^*(v)=\left\{ \begin{array}{ll}
0\ \ &\hbox{if} \ \ v \equiv 0,\\
\infty  \ \ &\hbox{if}  \ \ v \not \equiv 0, 
\end{array} \right.
\end{equation*} 
where $W^{-1,\infty}(\Omega)$ is the dual of $W^{1,1}_0(\Omega)$. The following lemma provides a formula for $E^*$.

\begin{lemma} \label{f^*lem} Let  $E$ be defined as in equation (\ref{FandG}). Then 

\begin{equation}\label{f^*}
E^*(V)=\left\{ \begin{array}{ll}
-\langle D f, V\rangle\ \ &\hbox{if} \ \ \varphi^0(x,V(x))\leq 1 \ \  \hbox{in} \ \ \Omega,\\
\infty  \ \ & \hbox{otherwise}.
\end{array} \right.
\end{equation}
\end{lemma}
{\bf Proof.} Suppose   
\begin{equation}\label{bigger}
\varphi^0(x,V(x))> 1
\end{equation}

on  a set $\omega \subset \Omega$ with positive Lebesgue measure. It follows from Lusin's theorem that for every $\epsilon>0$ there exists a compact set $Q \subset \omega$ such that $\mu(\omega \setminus Q)<\frac{\mu(\omega)}{2} $ and $V=\tilde{V}$ on $Q$, for some continuous function $\tilde{V} :\R^n \rightarrow \R^n$, where $\mu$ denotes the Lebesgue measure. In particular $\varphi^0(x,\tilde{V}(x))> 1$ for all $x\in Q$. Hence it follows from the definition of $\varphi^0$ that 
\[\forall x\in Q \ \ \ \ \exists P(x) \in S^{n-2} \ \ \hbox{such that}\ \ \varphi(x,P(x))< \tilde{V}(x)\cdot P(x).\]
Since $\tilde{V}$ and $\varphi$ are continuous, for every $x\in Q$ there exists $\epsilon_x$ such that 
\[\forall y\in B_{\epsilon_x}(x) \ \ \hbox{such that}\ \ \varphi(y, P(x)) < \tilde{V}(y)\cdot P(x). \]
Notice that
\[\{B_{\epsilon_x}(x): \ \ x\in Q\}\]
is an open cover for the compact set $Q$. Thus  there exists $z\in Q$ such that 
\[\mu (B_{\epsilon_{z}}(z)\cap Q)>0. \]
Now define $\bar{P} \in (L^{1}(\Omega))^n$ as follows 
\begin{equation*}\bar{P}=\left\{ \begin{array}{ll}
P(z) \ \ &\hbox{if} \ \ x \in B_{\epsilon_z}(z)\cap Q,\\
0  \ \ &\ \   \hbox{otherwise}.
\end{array} \right.
\end{equation*} 
Then we have  
\begin{eqnarray*}
E^*(V)&=&\sup_{P \in (L^{1}(\Omega))^n}  \left( \langle P, V \rangle-\int_{\Omega}\varphi(x, P+\nabla f)dx \right) \\
&=& -\langle \nabla f, V\rangle +\sup_{P \in (L^{1}(\Omega))^n} \left( \langle P, V \rangle-\int_{\Omega} \varphi(x, P)dx\right)\\
&\geq & -\langle \nabla f, V\rangle +\sup_{\lambda \in \R} \lambda \left( \langle \bar{P}, V \rangle-\int_{\Omega} \varphi(x, \bar{P})dx\right)\\
&=& -\langle \nabla f, V\rangle +\sup_{\lambda \in \R} \lambda  \int_{B_{\epsilon_z}(z)\cap Q}\left( V(x)\cdot P(z)- \varphi(x, P(z))\right) dx\\
&=& \infty.\\
\end{eqnarray*}
On the other hand if 
\[ \varphi^0(x,V(x))\leq 1,  \ \ \hbox{in}\ \ \Omega,\]
then 
\begin{equation}\label{NewLowerBound}
\varphi (x,P) \geq V(x) \cdot P  \ \  \ \ \forall \  x\in \Omega  \ \ \hbox{and}  \  \ P \in \R^n.
\end{equation}
Consequently

\begin{eqnarray*}
E^*(V)&=&\sup_{P \in (L^{1}(\Omega))^n}  \left( \langle P, V \rangle-\int_{\Omega}\varphi(x, P+\nabla f) dx\right) \\
&=& -\langle \nabla  f, V\rangle +\sup_{P \in (L^{1}(\Omega))^n} \left( \langle P, V \rangle-\int_{\Omega} \varphi(x, P)dx\right)\\
&=& -\langle \nabla f, V\rangle +\sup_{P \in (L^{1}(\Omega))^n}  \int_{\Omega}\left( V(x)\cdot P- \varphi(x, P)\right) dx\\
&=&  -\langle \nabla f, V\rangle .\\
\end{eqnarray*}
The proof is now complete. 
\hfill $\Box$

\vspace{.2cm}

Let $\nu_{\Omega}$ denote the outer unit normal vector to $\partial \Omega$, then for every $V\in (L^{\infty}(\Omega))^n$ with div$(V) \in L^n(\Omega)$ there exists a unique function $[V,  \nu_{\Omega}] \in L^{\infty}_{\mathcal{H}^{n-1}}(\partial \Omega)$ such that
\begin{equation}
\int_{\partial \Omega} [V,  \nu_{\Omega}] u d\mathcal{H}^{n-1}=\int_{\Omega} u \nabla \cdot V dx +\int_{\Omega}  V \cdot \nabla u dx, \ \ \forall u \in C^1(\bar{\Omega}).
\end{equation}

Moreover, for $u\in BV(\Omega)$ and $V\in (L^{\infty}(\Omega))^n$ with div$(V) \in L^n(\Omega)$, the linear functional $u\mapsto(V \cdot Du)$ gives rise to a Radon measure on $\Omega$, and
\begin{equation}\label{IBP0}
\int_{\partial \Omega} [V, \nu_{\Omega}]u d\mathcal{H}^{n-1}=\int_{\Omega} u \nabla \cdot V dx +\int_{\Omega}  (V \cdot D u), \ \ \forall u \in BV(\Omega),
\end{equation}
see \cite{Al, An} for a proof. See also Appendix C in \cite{And4} for a more recent exposition. 
Now define

\[\mathcal{V}:=\{V\in (L^{\infty}(\Omega))^n,\ \ \nabla \cdot V\equiv 0 \ \ \hbox{and }\varphi^0(x,V(x))\leq 1 \ \ \hbox{in} \ \ \Omega \}.\]
It follows from Lemma \ref{f^*lem} that the dual problem can be explicitly written as
\begin{equation*}\label{expDual} 
(P^*) \ \ \ \ \ \ \ \ \sup_{V \in \mathcal{V}}  \int _{\partial \Omega}f [V, \nu_{\Omega}] ds.  
\end{equation*}
where $\nu_{\Omega}$ is outward pointing unit normal vector on $\partial \Omega$.  The primal problem (P) may not have a solution, but the dual problem ($P^*$) always has a solution. This is a direct consequence of Theorem III.4.1 in \cite{ET}. Indeed it easily follows from \eqref{varphi.def1} that $I(v)=\int_{\Omega}\varphi (x, Dv)$ is convex, and $J: L^1(\Omega)\rightarrow \R$ with $J(p)=\int_{\Omega} \varphi(x, p)dx$ is continuous at $p=0$ (a consequence of $C_2$). Therefore the condition (4.8) in the statement of Theorem III.4.1 in \cite{ET} is satisfied.

\begin{proposition}\label{dualityGap}
Let $\Omega \subset \R^n$ be a bounded open set with Lipschitz boundary and assume $\varphi: \Omega \times \R^n \rightarrow \R$ be a continuous function satisfying the condition $(C_1)$ and $(C_2)$, and $f\in L^1(\partial \Omega)$.  Then there exists a divergence free vector field $T \in (L^{\infty}(\Omega))^n$ with $\varphi^0(x,T)\leq 1$ a.e. in $\Omega$ such that
\begin{eqnarray*}
\inf_{u\in W_f^{1,1}(\Omega)} \int_{\Omega} \varphi(x, Du) = \max_{V \in \mathcal{V}} \int _{\partial \Omega}f  [V, \nu_{\Omega}] ds= \int _{\partial \Omega}f  [T, \nu_{\Omega}] ds. 
\end{eqnarray*}
In particular, the dual problem $P^*$ has a solution $T \in \mathcal{V}$. 
\end{proposition}

In the above proposition $W^{1,1}_f$ denotes the space of functions in $W^{1,1}(\Omega)$ with trace $f$ on $\partial \Omega$. Proposition \ref{SameValue} follows directly from the following result. 
\begin{proposition} \label{MainProp}
Let $\Omega \subset \R^n$ be a bounded open set with Lipschitz boundary and assume $\varphi: \Omega \times \R^n \rightarrow \R$ be a continuous function satisfying the condition $(C_1)$ and $(C_2)$, and $f\in L^1(\partial \Omega)$.  Then 
\begin{eqnarray}\label{New}
\min _{u\in \mathcal{A}_f}  \left( \int_{\Omega} \varphi(x, Du)+\int_{\partial \Omega}\varphi(x,\nu_{\Omega})|f-u|ds \right)=   \inf_{u\in BV_f(\Omega)} \int_{\Omega} \varphi(x, Du) = \int _{\partial \Omega}f [T, \nu_{\Omega}] ds, 
\end{eqnarray}
where $T\in \mathcal{V}$ is a solution of the dual problem $P^*$ guaranteed by Proposition \ref{dualityGap}. 
\end{proposition}
{\bf Proof.} Let $u\in \mathcal{A}_f$ be a minimizer of \eqref{LTVProbRelaxed} and $T\in \mathcal{V}$ be a solution of the dual problem ($P^*$). Then  
\begin{eqnarray}
\int_{\Omega} \varphi(x, Du)&=& \int_{\Omega}\varphi(x,\frac{Du}{|Du|})|Du|\geq \int_{\Omega}T\cdot \frac{ Du}{|Du|}|Du|\\
&=& \int_{\Omega} T\cdot Du=\int_{\partial \Omega}u [T,\nu_{\Omega}] ds.\nonumber 
\end{eqnarray}
Now we conclude from \eqref{NewLowerBound} and the above inequality that 
\begin{eqnarray}\label{New11}
\int_{\Omega} \varphi(x, Du)+\int_{\partial \Omega}\varphi(x,\nu_{\Omega})|f-u|ds &\geq&  \int_{\partial \Omega}u [T,\nu_{\Omega}] ds +\int_{\partial \Omega} (f-u)[T,\nu_{\Omega}]ds \\
&=& \int_{\partial \Omega} f [T,\nu_{\Omega}]ds, 
\end{eqnarray}
and consequently 
\[\int_{\Omega} \varphi(x, Du)+\int_{\partial \Omega}\varphi(x,\nu_{\Omega})|f-u|ds  \geq \max_{V \in \mathcal{V}} \int _{\partial \Omega}f [V, \nu_{\Omega}] ds.\]
By Proposition \ref{dualityGap}, the above inequality holds also in the opposite direction, and hence \eqref{New} holds.  \hfill $\Box$

\vspace{.5 cm}

{\bf Proof of Theorem \ref{Structue}:}  Let $u\in \mathcal{A}_f$ be a minimizer of \eqref{LTVProbRelaxed} and $T\in \mathcal{V}$ be a solution of the dual problem ($P^*$). Then it follows from Proposition \eqref{MainProp} that 
\begin{eqnarray}\label{WeakIBP}
\int_{\partial \Omega}u [T,\nu_{\Omega}] ds=\int_{\Omega} \varphi(x, Du)&=& \int_{\Omega}\varphi(x,\frac{Du}{|Du|})|Du|\\
&\geq &  \int_{\Omega}T\cdot \frac{ Du}{|Du|}|Du| \nonumber \\
&=& \int_{\Omega} T\cdot Du=\int_{\partial \Omega}u [T,\nu_{\Omega}] ds,\nonumber 
\end{eqnarray}
and \eqref{directionInside} follows. On the other hand from Proposition \ref{dualityGap} and \eqref{WeakIBP} we conclude that 
\begin{eqnarray*}
\int_{\partial \Omega}f [T,\nu_{\Omega}]ds &\geq & \int_{\Omega} \varphi(x, Du)+\int_{\partial \Omega}\varphi(x,\nu_{\Omega})|f-u|ds \\
&\geq & \int_{\partial \Omega}u [T,\nu_{\Omega}] ds+ \int_{\partial \Omega}\varphi(x,\nu_{\Omega})|f-u|ds.
\end{eqnarray*}
Thus 
\[\int_{\partial \Omega}\varphi(x,\nu_{\Omega})|u-f|ds \leq  \int_{\partial \Omega}(f-u) [T,\nu_{\Omega}] ds. \]
Since $\varphi^0(x, T) \leq 1$ a.e. in $\Omega$, we have $\varphi(x,\nu_{\Omega}) \geq [T, \nu_{\Omega}]$. Hence \eqref{Boundary1} follows from the above inequality. \\

{\bf Acknowledgments.} I would like to thank Professors Robert L. Jerrard and Adrian Nachman for many insightful conversations. I would also like to thank the anonymous referee for their very careful reading of the paper and many useful comments which enormously improved the presentation of the paper. This research is partially supported by an start-up grant from University of California at Riverside.

\end{document}